# Identificación de Parámetros de Baterías de Ión-Litio: Un Estudio Comparativo de Diversos Modelos y Técnicas de Optimización para el Modelado de Baterías

*Parameter identification of lithium-ion batteries: a comparative study of various models and optimization techniques for battery modeling*

Johan Sebastian Suarez Sepúlveda

Mechatronics Engineer, Universidad Nacional de Colombia, Facultad de ingeniería, Bogotá, Colombia, jssuarezse@unal.edu.co , https://orcid.org/0009-0005-7898-7719. * Cra 45#56-07 apto 402, 3108027133

Edgar Hernando Sepúlveda-Oviedo

Ph.D. in Engineering, Laboratory for Analysis and Architecture of Systems – French National Centre for Scientific Research, Energy Management Systems Integration (ISGE), Toulouse, France. PROMES-CNRS (UPR-8521), Rambla de la thermodynamique, Tecnosud, 66100 Perpignan, France. Université de Perpignan Via Domitia, 52 Avenue Paul Alduy, 66860 Perpignan, France. ehsepulved@laas.fr, https://orcid.org/0000-0002-9947-3259

Bruno Jammes

Ph.D. in Engineering, Laboratory for Analysis and Architecture of Systems -, French National Centre for Scientific Research, Energy Management Systems Integration (ISGE), Toulouse, France. Université de Toulouse (UT), Toulouse, France. jammes@laas.fr, https://orcid.org/0009-0000-7793-8521

Corinne Alonso

Ph.D. in Engineering, Laboratory for Analysis and Architecture of Systems -, French National Centre for Scientific Research, Energy Management Systems Integration (ISGE), Toulouse, France. Université de Toulouse (UT), Toulouse, France. alonsoc@laas.fr, https://orcid.org/0000-0003-2547-255X

**Resumen**

En este trabajo se lleva a cabo un estudio comparativo de técnicas de optimización para la identificación de parámetros en modelos eléctricos equivalentes de baterías de ion-litio. Se emplea el modelo 2RC sobre un conjunto de doce baterías, utilizando cuatro bases de datos públicas provenientes de centros de investigación reconocidos. La metodología se estructura en cuatro etapas:

en primer lugar, se selecciona el modelo 2RC por su equilibrio entre precisión y simplicidad computacional; en segundo lugar, se recopilan datos experimentales de ciclos de carga y descarga; en tercer lugar, se aplican diversos métodos de optimización con el objetivo de minimizar el error entre los datos experimentales y los resultados estimados por el modelo; y finalmente, se evalúan la precisión, mediante el error cuadrático medio, y la eficiencia computacional, mediante el tiempo de ejecución. Se consideran métodos tradicionales, metaheurísticos y bioinspirados, entre ellos la optimización por mínimos cuadrados, el algoritmo de enjambre de partículas, el recocido simulado y diversas variantes inspiradas en procesos naturales. Se evidencia que las técnicas bioinspiradas permiten alcanzar una mayor precisión que los métodos tradicionales, sin un aumento significativo en el costo computacional. En particular, la optimización por enjambre de partículas muestra un desempeño superior en cuanto a exactitud y robustez frente a mínimos locales. Se concluye que la incorporación de estrategias de optimización avanzadas mejora significativamente la fidelidad de los modelos eléctricos equivalentes, lo que resulta fundamental para una estimación más precisa del estado de carga, el envejecimiento y la vida útil de baterías en aplicaciones críticas, tales como vehículos eléctricos y sistemas aeroespaciales.

**Palabras clave:** identificación de parámetros, baterías de ion-litio, modelos eléctricos equivalentes, técnicas de optimización, algoritmos bioinspirados

## Abstract (negrita, Calibri Light 11)

This work presents a comparative study of optimization techniques for parameter identification in equivalent electrical models of lithium-ion batteries. The 2RC model is applied to a set of twelve batteries using four publicly available datasets obtained from well-established research institutions. The methodology is structured in four stages: first, the 2RC model is selected due to its balance between physical interpretability and computational simplicity; second, experimental charge-discharge cycle data are collected; third, various optimization techniques are applied with the aim of minimizing the error between the experimental data and the response estimated by the model; and finally, accuracy is evaluated using the mean squared error, while computational efficiency is assessed through execution time. Traditional, metaheuristic, and bio-inspired optimization methods are considered, including least squares optimization, particle swarm optimization, simulated annealing, and several nature-inspired variants. It is demonstrated that bio-inspired techniques achieve greater accuracy than traditional methods, without a significant increase in computational cost. In particular, particle swarm optimization shows superior performance in terms of precision and robustness against local minima. It is concluded that the integration of advanced optimization strategies significantly enhances the fidelity of equivalent electrical models, which is essential for accurate estimation of internal states such as state of charge, aging, and service life in lithium-ion batteries used in electric vehicles and aerospace systems.

**Keywords:** parameter identification, lithium-ion batteries, equivalent electrical models, optimization techniques, bio-inspired algorithms

# 1. INTRODUCCIÓN

Lithium-ion batteries are considered a cornerstone technology for modern energy storage applications due to their high energy and power density, long cycle life, and reduced self-discharge rates [22]. These characteristics have made them the dominant energy source in electric vehicles, portable electronics, and aerospace systems. Their integration is also essential to support the large-scale deployment of renewable energy sources, where energy storage systems play a key role in stabilizing power grids [27]. As the demand for efficient and safe battery systems increases, significant attention has been directed toward the accurate modeling of lithium-ion battery behavior, particularly under variable charge and discharge conditions.

To ensure effective battery management and extend operational life, internal states such as the state of charge (SOC), state of health (SOH), and thermal behavior must be estimated in real time with high accuracy. The battery management system relies on mathematical models to perform this task. Among the existing modeling strategies, equivalent circuit models (ECMs) are widely used due to their balance between simplicity and performance [2]. ECMs are composed of passive electrical components—resistors, capacitors, and voltage sources—that replicate the dynamic behavior of batteries during operation. Compared to electrochemical models, ECMs require fewer parameters, are computationally lighter, and are easier to implement in embedded systems [1].

Despite their advantages, the performance of ECMs is strongly dependent on the accurate identification of their parameters. These parameters are not constant; they vary with temperature, current rates, aging, and usage history [11]. The identification process is challenging because it involves solving nonlinear optimization problems in high-dimensional search spaces. Traditional methods such as least squares optimization are widely used due to their speed and simplicity; however, they often get trapped in local minima and fail to capture the complexity of battery behavior under dynamic conditions [25].

To overcome these limitations, various metaheuristic and bio-inspired optimization techniques have been explored in the literature. Algorithms such as Particle Swarm Optimization, Genetic Algorithms, Simulated Annealing, and more recent nature-inspired methods like the Marine Predators Algorithm [6] have demonstrated superior performance in solving complex parameter estimation problems [29]. These methods offer a better exploration of the solution space, making them particularly suitable for nonlinear systems such as lithium-ion batteries. Nonetheless, these approaches may incur higher computational costs, which is a critical consideration in real-time applications.

In response to these challenges, this study contributes to the field by conducting a comparative analysis of several optimization techniques for parameter identification in lithium-ion battery models. The analysis focuses on the 2RC equivalent circuit model, which balances physical interpretability and computational simplicity. Publicly available datasets from multiple well-established research centers are employed to ensure the robustness and generality of the findings. Unlike most studies available in the state of the art, this work emphasizes the use of limited discharge cycle data, with the specific aim of enabling fast parameter identification suitable for embedded systems. This approach is particularly relevant for applications where continuous measurement of the full discharge cycle is impractical due to time, energy, or resource constraints. The main contributions of this work are as follows:

    I.     a systematic evaluation of traditional, metaheuristic, and bio-inspired optimization algorithms applied to battery parameter estimation;

    II.    a performance comparison based on estimation accuracy and computational cost using standardized metrics such as the mean squared error and execution time;

    III.   the identification of techniques that provide an optimal trade-off between accuracy and efficiency for real-time applications in energy storage systems; and

    IV.   the introduction of a data-efficient modeling strategy that supports the development of lightweight estimation schemes for battery management in embedded platforms.

This study is structured as follows. Section 2 provides a brief overview of common battery modeling techniques to motivate the selection of the 2RC equivalent circuit model. Section 3 describes the methodology used, including the model configuration, datasets, and implementation environment. Section 4 presents the results obtained from the optimization techniques. Section 5 discusses the findings in terms of accuracy and computational cost. Finally, Section 6 summarizes the conclusions and outlines directions for future work.

## 2. RELATED WORK

The development of accurate battery models and robust parameter identification methods has received increasing attention due to their importance for real-time battery management in electric vehicles, renewable energy storage, and aerospace systems. This section reviews several recent works focused on parameter estimation of lithium-ion batteries using optimization algorithms, highlighting the key differences and contributions of our study.

In [23], Tian et al. proposed a Dragonfly Algorithm (DA)-based identification scheme using a first-order RC model. While the approach achieved high accuracy, it was limited by the simplicity of the model and did not compare favorably against more advanced circuit representations like the 2RC model. Our study advances this by using a dual RC configuration and demonstrating better performance using PSO, which allows for more detailed modeling of dynamic behaviors such as charge transfer and diffusion.

Ferahtia et al. [7] explored parameter identification using the Artificial Ecosystem Optimization (AEO) algorithm on a Shepherd model. Their results show high identification efficiency; however, the Shepherd model is empirical and less suited for integration in embedded systems due to its limited physical interpretability. In contrast, our work adopts the 2RC equivalent circuit model, which provides a better balance between simplicity and fidelity, and is more suitable for real-time embedded implementations.

A recent work by Wang et al. [26] introduced a classification model-assisted Bayesian optimization (CMABO) approach to overcome convergence failures in simulations with P2D models. Although this framework significantly reduces simulation failures and improves identification efficiency, it requires the use of computationally intensive electrochemical models. Our methodology, instead, prioritizes data-efficiency and real-time applicability by using a simpler but highly expressive 2RC model, optimized through PSO with significantly lower computational cost.

Another notable contribution is from Gu et al. [8], who proposed a sensitivity-oriented stepwise optimization (SSO) method for electrochemical parameter identification. Although their method yields high precision for degradation-related parameters, it depends on extensive experimental setups and sequential testing. Our approach eliminates the need for such complex testing protocols, achieving competitive results using public datasets and standard discharge profiles.

Finally, Li et al. [14] performed a large-scale benchmark of 78 metaheuristic methods for identifying electrochemical model parameters. Their findings revealed strong results for teaching-learning-based optimization, but at the expense of moderate to high computational burden. Our results confirm that PSO can achieve similar accuracy levels with far less computational complexity, making it more suitable for lightweight applications.

Unlike these studies, our work stands out in several aspects: (1) it adopts a second-order ECM for better dynamic representation of lithium-ion batteries; (2) it performs a comprehensive comparative evaluation of traditional, metaheuristic, and bio-inspired methods; (3) it identifies PSO as a method offering optimal balance between accuracy and execution time; and (4) it demonstrates the feasibility of achieving sub-millisecond errors using only partial discharge data, a crucial feature for embedded system applications.

## 3. METHODOLOGY

This section outlines the methodological approach used for parameter identification in lithium-ion battery models. The methodology is divided into three main components. First, we justify the selection of the 2RC equivalent circuit model based on its balance between accuracy and computational efficiency, making it suitable for embedded applications. Second, we describe the comparative framework for optimization techniques, including both traditional and bio-inspired methods, selected for their potential to improve parameter estimation performance under nonlinear and high-dimensional conditions. Finally, we present the experimental setup, including dataset sources, implementation tools, and evaluation metrics.

### 3.1. Model Selection and Justification

Battery modeling methods can be broadly classified into four categories: empirical models, electrochemical models, data-driven models, and equivalent circuit models (ECMs). Empirical models simplify battery behavior using curve-fitting techniques and low-order polynomials. While easy to implement, they often suffer from poor generalization due to measurement noise and nonlinear effects such as hysteresis and memory phenomena [18]. Electrochemical models, such as the pseudo two-dimensional (P2D) model or the Single Particle Model (SPM), provide detailed insight into internal physicochemical processes. However, they require the solution of partial differential equations and estimation of a large number of parameters, which makes them unsuitable for real-time and embedded applications [4].

Data-driven models—including neural networks, support vector machines, and fuzzy logic systems—learn mappings from input data to internal states without relying on physical principles. These approaches have demonstrated good performance in state estimation and fault detection but tend to lack physical interpretability and require large datasets for training. Their practical implementation in embedded systems is still limited by memory and processing constraints [12].

In contrast, equivalent circuit models provide a compromise between accuracy and computational efficiency. ECMs approximate battery dynamics using passive electrical elements such as resistors, capacitors, and voltage sources. These models are able to capture critical phenomena such as polarization and diffusion while maintaining simplicity for embedded deployment [4]. Several topologies have been proposed, including Rint, Thevenin (1RC), PNGV, and 2RC models. Among them, the 2RC model is particularly effective in reproducing both fast and slow dynamic responses, making it well suited for robust state estimation under varying operating conditions [18].

For these reasons, this study focuses exclusively on the 2RC equivalent circuit model. Unlike more complex models, it allows for accurate parameter identification using only a limited portion of the discharge cycle. This is particularly advantageous for real-time applications in embedded battery management systems where complete data collection is often impractical.

### *2RC Model or Dual Polarization Model*

The 2RC model (See Fig. 1.) is widely used for online estimation of the state of charge due to its excellent trade-off between accuracy and computational complexity. It is particularly effective in capturing the dynamic polarization effects present in lithium-ion batteries. This model includes two parallel resistor-capacitor branches that represent both fast and slow transient processes associated with charge transfer and concentration gradients within the cell [15].

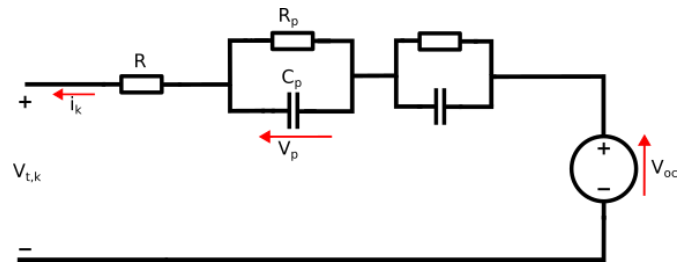

Figure 1. 2RC Model or Dual Polarization Model

Compared to simpler models such as the Thevenin (1RC) [9] or Rint configurations [28], the 2RC model demonstrates superior fidelity in representing internal battery dynamics under dynamic load profiles. Recent enhancements, such as those introduced in modified PNGV variants, aim to reduce estimation errors caused by state of charge variability and temperature dependence. For this reason, the 2RC model is considered one of the most reliable choices for embedded systems requiring real-time performance and accurate modeling [17].

The equivalent circuit consists of an internal resistance $R$, two $RC$ branches $(R_1, C_1)$ and $(R_2, C_2)$, and an open-circuit voltage source $V_{OC}$. These components together replicate the battery voltage behavior under varying load conditions. The governing equation is (1).

$$V_{t,k} = V_{OC}(SOC_k) - V_{1,k} - V_{2,k} - I_k \cdot R \qquad (1)$$

Where $V_{OC}(SOC_k)$ is the open-circuit voltage at sample $k$, $V_{1,k}$ and $V_{2,k}$ are the voltage drops across the RC branches, $I_k$ is the terminal current, and $R$ is the internal resistance. The values of these parameters depend strongly on the state of charge and the operating temperature of the battery, making precise parameter identification a key requirement for model performance [17].

### 3.2. Optimization Strategy Comparison

Accurate parameter identification is essential to the performance of equivalent circuit models (ECMs), as it directly influences their ability to reproduce lithium-ion battery behavior under dynamic and nonlinear operating conditions. Given the highly nonlinear, multidimensional nature of the parameter space, selecting an appropriate optimization strategy is a critical step in battery modeling.

Historically, traditional deterministic methods such as least squares optimization have been widely used due to their simplicity and low computational cost [24]. However, these methods often suffer from limited robustness when dealing with noisy or incomplete data, and can converge to local minima. To address these limitations, numerous metaheuristic and bio-inspired algorithms have been proposed that offer improved global search capabilities and adaptability to high-dimensional spaces.

Among these, Particle Swarm Optimization (PSO) has gained prominence due to its ability to balance exploration and exploitation efficiently. PSO has been successfully applied to estimate parameters in ECMs, achieving better accuracy and faster convergence than deterministic techniques under various dynamic load profiles [5]. Comparative studies have demonstrated that PSO outperforms traditional methods in both estimation error and convergence speed [16].

More recent innovations include nature-inspired metaheuristics such as the Marine Predators Algorithm (MPA), which mimics the intelligent foraging behavior of marine predators and has shown strong results in complex energy optimization problems [6], as well as its application to photovoltaic systems [17]. While MPA shows promise in battery modeling, it often involves more computational overhead than simpler approaches like PSO.

In this study, a comparative evaluation of optimization methods was conducted, including Least Squares (LS), Marine Predators Algorithm (MPA), and Particle Swarm Optimization (PSO). Each method was tested on identical datasets and evaluated using two metrics: execution time (ET) and mean squared error (MSE) between the estimated and actual battery voltage curves. MSE is widely accepted as a standard accuracy metric for model validation in regression problems [13].

Based on the results, PSO was selected as the primary optimization method in this work due to its superior trade-off between accuracy and computational efficiency. It demonstrated a lower MSE than LS and comparable performance to MPA with significantly reduced computation time—making it ideal for embedded and real-time applications. These findings align with prior studies emphasizing the reliability of PSO for battery model calibration tasks [19].

Given its central role in our implementation, a more detailed overview of the PSO algorithm is provided below, including its operational principles and key advantages in the context of battery parameter estimation.

*Particle Swarm Optimization*

Particle Swarm Optimization (PSO) is a population-based metaheuristic inspired by the collective behavior observed in natural systems such as bird flocking and fish schooling. Originally introduced by Kennedy and Eberhart (1995), PSO simulates a group of particles—each representing a candidate solution—that move through the search space by updating their velocity and position based on both individual experience and collective knowledge of the swarm. The position of each particle is adjusted iteratively according to the best solution found by itself (personal best) and by the entire swarm (global best), allowing for a balance between exploration of new areas and exploitation of promising regions.

In the context of parameter identification for lithium-ion battery models, PSO offers several critical advantages. First, it does not require gradient information or specific model assumptions, making it highly suitable for non-convex and multidimensional optimization problems such as those encountered in equivalent circuit models (ECMs). Second, PSO is computationally efficient, with a relatively simple algorithmic structure that makes it well suited for embedded systems with limited resources. Third, PSO has demonstrated strong convergence properties in diverse engineering applications, including battery modeling, due to its ability to avoid premature convergence and local minima traps [13].

In recent literature, PSO has consistently outperformed traditional methods like Least Squares (LS) in terms of estimation accuracy and robustness under dynamic operating conditions [6]. Comparative studies have also shown that, while more advanced metaheuristics such as the Marine Predators Algorithm (MPA) can provide slightly better accuracy, they often come at the cost of increased computational complexity and longer convergence times [29].

In this study, PSO was implemented as the principal optimization strategy for estimating the parameters of the 2RC model. The objective function minimized by PSO is the mean squared error (MSE) between the model-predicted voltage and the actual experimental data. Key parameters of the PSO algorithm—including inertia weight, cognitive coefficient, and social coefficient—were tuned empirically based on preliminary simulations to ensure stable convergence and reduced execution time. As shown in Section 4, PSO achieved a favorable balance between estimation accuracy and computational efficiency, outperforming LS and matching or exceeding the performance of MPA in most cases. These characteristics make PSO an optimal choice for real-time applications such as onboard battery management systems, where both speed and reliability are essential.

## 3.3. Experimental setup

This section details the experimental framework used to evaluate and compare the performance of optimization strategies for parameter identification in the 2RC equivalent circuit model.

*Dataset sources*

To ensure the robustness and generality of the analysis, datasets were sourced from four publicly available and widely referenced repositories: NASA, CALCE, Oxford, and HNEI. These datasets include full discharge cycles for lithium-ion cells under various load conditions and temperatures, providing a representative basis for testing parameter estimation techniques. The datasets were originally compiled and harmonized in [3], and are accessible via the Battery Archive repository. In total, discharge data from 12 different cells were used and cycled until the end of their lifespan, resulting in 1,155 identifications. These identifications were used to compare the performance of the eight optimization methods.

*Implementation tools*

All simulations and model implementations were carried out in MATLAB, leveraging its Optimization and Global Optimization Toolboxes. The equivalent circuit models and the parameter estimation routines were implemented using vectorized functions to ensure efficient computation. For reproducibility and to support real-time feasibility analysis, the experiments were run on a machine with standard embedded-system-comparable resources: an Intel i7 processor and 16GB of RAM. All optimization algorithms used in the comparative evaluation—Least Squares (LS), Particle Swarm Optimization (PSO), and others—were coded consistently under a unified framework to ensure fair comparison.

*Evaluation metrics*

The performance of each optimization algorithm was assessed using two complementary metrics: mean squared error (MSE) and execution time (ET).

- **MSE** quantifies the average of the squared differences between the predicted voltage from the model and the experimental measurements. It is widely adopted for model validation in system identification and regression contexts [13].

- **ET** corresponds to the total time, in seconds, taken by the optimization algorithm to converge to a final solution. This metric is particularly relevant for real-time and embedded applications where computational constraints are strict.

These metrics were selected to reflect both the estimation accuracy and computational feasibility of each optimization method under study.

## 4. RESULTS

Eight optimization techniques were evaluated for the purpose of accurately identifying the parameters of the 2RC equivalent circuit model for lithium-ion batteries. The comparison was based on two performance metrics: Mean Squared Error (MSE) as a measure of estimation accuracy, and Execution Time (ET) as an indicator of computational efficiency. The results are summarized in Table 1 and further illustrated through a comparative voltage-time plot in Figure 2.

**Table 1. Computational Cost and Accuracy of the Different Optimization Methods**

| Method | MSE | ET (s) |
|---|---|---|
| Least Squares | $1.698 \times 10e^{-6}$ | 0.15 |
| **Particle Swarm** | $\mathbf{3.577 \times 10e^{-7}}$ | **0.56** |
| Simulated Annealing | $5.9719 \times e^{-7}$ | 0.61 |
| Genetic Algorithm | $4.66 \times 10e^{-6}$ | 1.47 |
| Golf Field | $7.074 \times 10e^{-6}$ | 1.25 |
| Australian Dingo | $3.907 \times 10e^{-7}$ | 2.98 |
| Mexican Axolotl | $1.23 \times 10e^{-6}$ | 2.91 |
| Spider Jumping | $7.234 \times 10e^{-5}$ | 4.2 |

Fuente: elaboración propia.

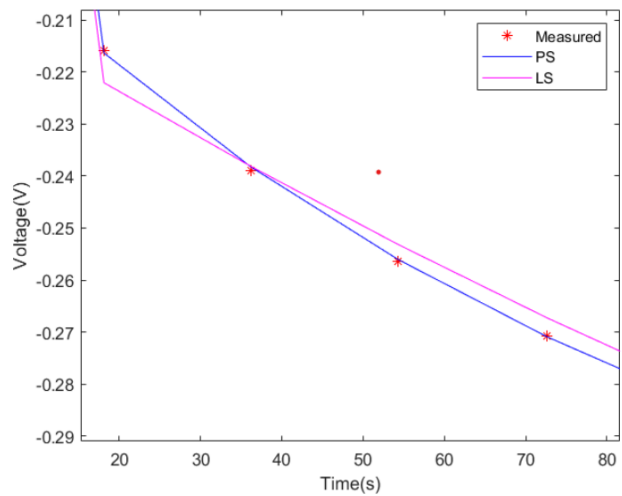

Figure 2. Comparative of the Accuracy of the Optimization Methods Least Squares (LS) and Particle Swarm Optimization (PS).

Among all tested methods, Particle Swarm Optimization (PSO) demonstrated the best overall performance. It achieved the lowest MSE at $3.577 \times 10e^{-7}$, with an execution time of just 0.57 seconds in each identification. This represents a substantial improvement in accuracy over the Least Squares (LS) method, which yielded a higher MSE of $1.698 \times 10e^{-6}$, albeit with the shortest execution time (0.15s). Other methods such as Simulated Annealing and Genetic Algorithm offered modest improvements but were not able to outperform PSO in either accuracy or speed.

In contrast, more complex bio-inspired methods such as Mexican Axolotl and Spider Jumping exhibited significantly longer execution times (2.91 s and 4.2 s, respectively) and much higher MSE values, making them impractical for real-time embedded applications where resource constraints are critical.

The accuracy advantage of PSO is further confirmed in Figure 2, which displays a voltage-time comparison between measured experimental values and model predictions using LS and PSO. It is evident that the voltage curve generated by PSO adheres more closely to the actual measurements, particularly in regions with sharp changes, underscoring PSO's superior parameter estimation capability under dynamic conditions.

In conclusion, Particle Swarm Optimization achieved the most favorable trade-off between precision and computational efficiency, making it the most suitable method for battery parameter identification in this study and a strong candidate for use in real-time applications such as embedded battery management systems.

## 5. DISCUSSION

The results presented in Section 4 demonstrate a clear advantage of Particle Swarm Optimization (PSO) over both traditional deterministic methods and other metaheuristic algorithms in the context of parameter identification for lithium-ion battery models. In this section, the implications of these findings are discussed in detail, including key observations, methodological insights, and comparison with related work.

### *Key Findings and Interpretation*

From the comparative analysis, PSO yielded the lowest mean squared error ($3.57 \times 10e^{-7}$) with a moderate execution time of 0.56 seconds. This performance surpasses traditional methods such as Least Squares (LS), which while computationally efficient (0.15 s), exhibited substantially higher error ($1.68 \times 10e^{-6}$). PSO also outperformed other advanced algorithms such as Simulated Annealing and Genetic Algorithms, both of which showed higher MSEs and similar or longer execution times.

Furthermore, the visualization in Figure 2 confirms the superior fitting capability of PSO across the voltage-time domain. Predictive curve of the PSO adheres more closely to the measured data, especially in regions with rapid voltage change, indicating better generalization under dynamic load conditions.

The effectiveness of PSO lies in its balance between global search (exploration) and local refinement (exploitation), which allows it to avoid premature convergence and efficiently handle the multidimensional optimization space inherent to battery modeling.

### Challenges

Despite its robust performance, PSO and other heuristic methods are not without challenges. One major difficulty is the sensitivity of the algorithm to parameter tuning—namely, the inertia weight and cognitive/social coefficients—which must be calibrated to avoid stagnation or erratic behavior. Additionally, while PSO performs well on average, its stochastic nature implies that occasional runs may yield suboptimal results without sufficient swarm diversity.

Another challenge is model-data compatibility. The 2RC model, while effective, abstracts many electrochemical phenomena. Thus, for high-fidelity applications such as health prognostics or thermal coupling, ECMs may need to be combined with other modeling strategies or augmented with adaptive parameter correction.

### Study Limitations

The design of this study introduces certain limitations that must be acknowledged:

1. **Restricted Model Class**: The analysis was limited to the 2RC equivalent circuit model. Although widely used, it does not capture all internal dynamics such as temperature effects or aging-induced changes in impedance.

2. **Dataset Scope**: The work was conducted using a curated subset of 15 battery discharge cycles from four publicly available datasets. While this ensures generalizability, it does not cover all chemistries or degradation states.

3. **Static Evaluation Criteria**: The evaluation was based on MSE and ET under a fixed load profile. Performance under varied real-world cycles, including charge events or hybrid profiles, remains to be assessed.

4. **Limited Hardware Benchmarking**: The computational evaluation was conducted on a desktop processor. For embedded applications, performance on low-power microcontrollers must be explicitly validated.

### Future Work

To build on the present findings, several future research directions are proposed:

- **Hybrid Modeling**: Combine ECMs with machine learning to adapt model parameters online using real-time sensor feedback. This would improve robustness against aging and thermal variations. [10]

- **Multiscale Optimization**: Integrate coarse-to-fine techniques such as PSO-LM hybrids [21] to further reduce convergence time while maintaining global accuracy.
- **Model Generalization**: Extend the current framework to incorporate variable load cycles, partial charges, and low-temperature operation to ensure broader applicability.
- **Embedded Implementation**: Validate the PSO-based estimation algorithm on real-time hardware such as STM32 or Raspberry Pi platforms to assess its practical feasibility for Battery Management Systems (BMS).
- **Uncertainty Quantification**: Introduce probabilistic frameworks or Bayesian optimization to quantify confidence intervals around parameter estimates, which is critical for safety-critical applications such as aerospace.
- **Scientific Landscape Analysis**: Perform a computational analysis of the research landscape using bibliometric mapping and topic modeling tools. This will help identify emerging trends, research gaps, and related subfields, thereby refining the scope and positioning of future contributions [20].

## Comparative Context with State-of-the-Art

Our findings align with and extend the results from several recent studies. For instance, Tian et al. [23] employed a Dragonfly Algorithm for parameter estimation but restricted their model to a simple first-order RC circuit, limiting dynamic representation. Ferahtia et al. [7] used the AEO algorithm on the Shepherd model, achieving high accuracy but with limited relevance for embedded deployment due to the empirical nature of the model.

Other works, such as Wang et al. [26], applied a Bayesian framework to accelerate convergence in electrochemical models. However, the computational cost remains prohibitive for lightweight systems. Gu et al. [8] proposed a sensitivity-oriented optimization scheme but required elaborate testing protocols for parameter decomposition, which are not always available in practice.

Compared to these studies, our work offers a unique combination of lightweight model complexity (2RC), open-source data use, real-time feasibility, and multi-algorithm comparison under identical conditions. These contributions support its adoption in embedded battery diagnostics, especially where full discharge cycles are not feasible to acquire.

## 6. CONCLUSIONES

This study conducted a comprehensive comparative analysis of parameter identification methods for lithium-ion battery models using various optimization techniques. By focusing on the 2RC equivalent circuit model and leveraging publicly available datasets from multiple research institutions, the work ensures both methodological rigor and practical relevance.

The results demonstrate that Particle Swarm Optimization (PSO) outperforms both traditional deterministic approaches, such as Least Squares (LS), and more computationally intensive bio-inspired algorithms. PSO achieved the lowest mean squared error ($3.57 \times 10^{-7}$) with a moderate execution time (0.56 s), making it a compelling choice for real-time applications in battery management systems. The superior fitting accuracy of PSO was evident in dynamic regions of the voltage profile, highlighting its robustness in handling nonlinear and time-varying behavior.

The main contributions of this study include: (1) a systematic evaluation of eight optimization algorithms applied to the parameter identification task; (2) the identification of PSO as the most balanced method in terms of accuracy and computational cost; (3) the validation of the 2RC model as a reliable yet computationally efficient framework for lithium-ion battery modeling; and (4) the demonstration that accurate parameter estimation can be achieved using limited discharge data, thus enabling data-efficient and embedded-friendly modeling strategies.

These findings provide a solid foundation for the deployment of PSO-based estimation schemes in embedded platforms, and support future research aimed at hybrid modeling, uncertainty quantification, and real-time diagnostics for advanced battery systems.

## REFERENCIAS

[1] Allagui, A., Freeborn, T. J., Elwakil, A. S., Fouda, M. E., Maundy, B. J., Radwan, A. G., Said, Z, & Abdelkareem, M. A. (2018). Review of fractional-order electrical characterization of supercapacitors. *Journal of Power Sources, 400*, 457–467.

[2] Bandhauer, T. M., Garimella, S., & Fuller, T. F. (2011). A critical review of thermal issues in lithium-ion batteries. *Journal of The Electrochemical Society, 158*, R1.

[3] Battery Archive. (2025). *Battery Archive Dataset*. Retrieved from http://batteryarchive.org

[4] Bedwal, K., & Moulik, B. (2024). Review on equivalent circuit and data driven approaches for electric vehicle battery models. *2nd International Conference on Advancements and Key Challenges in Green Energy and Computing (AKGEC)*, 1–6.

[5] Dziechciaruk, G., Ufnalski, B., & Grzesiak, L. (2017). Parameter estimation for equivalent electrical model of lithium-ion cell. *2017 19th European Conference on Power Electronics and Applications (EPE'17 ECCE Europe)*, P.1–P.9.

[6] Faramarzi, A., Heidarinejad, M., Mirjalili, S., & Gandomi, A. H. (2020). Marine Predators Algorithm: A nature-inspired metaheuristic. *Expert Systems with Applications, 152*, 113377.

[7] Ferahtia, S., Djeroui, A., Rezk, H., Chouder, A., Houari, A., & Machmoum, M. (2021). Optimal parameter identification strategy applied to lithium-ion battery model. *International Journal of Energy Research, 45*, 16741–16753.

[8] Gu, Y., Wang, J., Chen, Y., Zheng, K., Deng, Z., & Chen, Q. (2022). Electrochemical parameter identification for lithium-ion battery sources in self-sustained transportation energy systems. *IEEE Transactions on Industry Applications, 60*, 1240–1254.

[9] Hannan, M. A., Lipu, M. S. H., Hussain, A., & Mohamed, A. (2017). A review of lithium-ion battery state of charge estimation and management system in electric vehicle applications: Challenges and recommendations. *Renewable and Sustainable Energy Reviews, 78*, 834–854.

[10] Jammes, B., Sepúlveda-Oviedo, E. H., & Alonso, C. (2025). Suivi temps réel du SoH de batteries lithium-ion. *Symposium de Génie Électrique (SGE 2025) *, 1–5.

[11] Jiang, Y., Xia, B., Zhao, X., Nguyen, T., Mi, C., & de Callafon, R. A. (2017). Data-based fractional differential models for non-linear dynamic modeling of a lithium-ion battery. *Energy, 135*, 171–181.

[12] Kok, S. X., Hau, L. C., Lim, Y. S., Wong, J., Chua, K. H., & Stella, M. (2024). Investigation of data-driven models for robust lithium-ion battery modeling. *IEEE Sustainable Power and Energy Conference (iSPEC)*, 1–6.

[13] Li, Y., Jiao, J., Yang, Y., & Ji, P. (2022). Parameter estimation of equivalent circuit model for lithium batteries. In Xie, Q., Zhao, L., Li, K., Yadav, A., & Wang, L. (Eds.), *Advances in Natural

Computation, Fuzzy Systems and Knowledge Discovery. Lecture Notes on Data Engineering and Communications Technologies* (Vol. 89, pp. xx–xx). Springer.

[14] Li, Y, Liu, G., Deng, W., & Li, Z. (2024). Comparative study on parameter identification of an electrochemical model for lithium-ion batteries via meta-heuristic methods. *Applied Energy.*.

[15] Liu, S., Deng, D., Wang, S., Luo, W., Takyi-Aninakwa, P., Qiao, J., Li, S., Jin, S., & Hu, C. (2023). Dynamic adaptive square-root unscented Kalman filter and rectangular window recursive least square method for the accurate state of charge estimation of lithium-ion batteries. *Journal of Energy Storage, 67*, 107603.

[16] Mohamed, M. A. A., Yu, T. F., Ramsden, G., Marco, J., & Grandjean, T. (2023). Advancements in parameter estimation techniques for 1RC and 2RC equivalent circuit models of lithium-ion batteries: A comprehensive review. *Journal of Energy Storage, 113*, 115581.

[17] Peng, J., Meng, J., Wu, J., Deng, Z., Lin, M., Mao, S., & Stroe, D. I. (2023). A comprehensive overview and comparison of parameter benchmark methods for lithium-ion battery application. *Journal of Energy Storage, 71*, 108197.

[18] Riemann, B. J. C., Li, J., Adewuyi, K., Landers, R. G., & Park, J. (2021). Control-oriented modeling of lithium-ion batteries. *Journal of Dynamic Systems, Measurement, and Control, 143*, 021002.

[19] Sangwan, V., Sharma, A., Kumar, R., & Rathore, A. K. (2016). Estimation of battery parameters of the equivalent circuit models using meta-heuristic techniques. *2016 IEEE 1st International Conference on Power Electronics, Intelligent Control and Energy Systems (ICPEICES)*, 1–6.

[20] Sepúlveda-Oviedo, E. H., Travé-Massuyès, L., Subias, A., Pavlov, M., & Alonso, C. (2023). Fault diagnosis of photovoltaic systems using artificial intelligence: A bibliometric approach. *Heliyon, 9*, e21491.

[21] Shen, W., & Li, H. (2017). Multi-scale parameter identification of lithium-ion battery electric models using a PSO-LM algorithm. *Energies, 10*, 1–18.

[22] Thakur, A. K., Prabakaran, R., Elkadeem, M., Sharshir, S. W., Arıcı, M., Wang, C., Zhao, W., Hwang, J. Y., & Saidur, R. (2020). A state of art review and future viewpoint on advance cooling techniques for Lithium–ion battery system of electric vehicles. *Journal of Energy Storage, 32*, 101771.

[23] Tian, J., Yin, X., Pan, T., Zhang, X., Yang, D., & Ni, L. (2024). Parameter identification of lithium-ion battery using Dragonfly Algorithm. *2024 IEEE 25th China Conference on System Simulation Technology and its Application (CCSSTA)*, 565–569.

[24] Tian, N., Wang, Y., Chen, J., & Fang, H. (2017). On parameter identification of an equivalent circuit model for lithium-ion batteries. *IEEE Conference on Control Technology and Applications (CCTA)*, 187–192.

[25] Waag, W., Käbitz, S., & Sauer, D. U. (2013). Application-specific parameterization of reduced order equivalent circuit battery models for improved accuracy at dynamic load. *Measurement, 46*, 4085–4093.

[26] Wang, B., He, Y., Liu, J., & Luo, B. (2023). Fast parameter identification of lithium-ion batteries via classification model-assisted Bayesian optimization. *Energy. *

[27] Wang, Y., Liu, B., Li, Q., Cartmell, S., Ferrara, S., Deng, Z. D., & Xiao, J. (2015). Lithium and lithium-ion batteries for applications in microelectronic devices: A review. *Journal of Power Sources, 286*, 330–345.

[28] Wu, M., Qin, L., Wu, G., Huang, Y., & Shi, C. (2021). State of charge estimation of power lithium-ion battery based on a variable forgetting factor adaptive Kalman filter. *Journal of Energy Storage, 41*, 102841.

[29] Yousri, D., Babu, T. S., Beshr, E., Eteiba, M. B., & Allam, D. (2020). A robust strategy based on marine predators' algorithm for large scale photovoltaic array reconfiguration to mitigate the partial shading effect on the performance of PV system. *IEEE Access, 8*, 112407–112426.

---

**Edgar Hernando Sepúlveda-Oviedo**

- **Academic Background:** Mechatronics Engineer from the National University of Colombia (2016), where he graduated at the top of his class. He completed his Master's degree in Industrial Automation at the same university (2016–2019), graduating with honors and receiving a distinguished thesis award. He obtained his Ph.D. from Université Paul Sabatier and the Laboratory for Analysis and Architecture of Systems at the French National Centre for Scientific Research (LAAS-CNRS) in Toulouse, France. He then worked as a postdoctoral researcher at LAAS-CNRS, within the "Intégration de systèmes de gestion de l'énergie" (ISGE) team, focusing on the application of artificial intelligence to embedded energy systems. He is currently an Associate Professor at the University of Perpignan and a member of the PROMES laboratory (CNRS), where he continues his research in intelligent diagnostic and control for next-generation solar energy systems.

- **Relevant Publications:**

  -Sepúlveda-Oviedo, E.H., Travé-Massuyès, L., Subias, A., Pavlov, M., Alonso, C. (2023). Fault diagnosis of photovoltaic systems using artificial intelligence: A bibliometric approach. *Heliyon*, 9, e21491.

  -Sepúlveda-Oviedo, E. H., Travé-Massuyès, L., Subias, A., Alonso, C., Pavlov, M. (2022). Feature extraction and health status prediction in PV systems. *Advanced Engineering Informatics*, 53, 101696.

  -Sepúlveda-Oviedo, E.H., Travé-Massuyès, L., Subias, A., Pavlov, M., Alonso, C. (2024). An ensemble learning framework for snail trail fault detection and diagnosis in photovoltaic modules. *Engineering Applications of Artificial Intelligence*, 137 (Part A), 1090688.

- **Institutional Affiliations:** National University of Colombia (Bogotá, Colombia): Assistant Professor in the Faculty of Engineering (2017). Industrial Automation Research Group (GAUNAL), National University of Colombia: Associate Researcher. Master's in Industrial Automation, National University of Colombia (2016–2019): Graduate Student. Laboratory for Analysis and Architecture of Systems (LAAS-CNRS, Toulouse, France): Doctoral

Researcher (2020–2023); Postdoctoral Researcher (2023–2024). Institut National des Sciences Appliquées de Toulouse (INSA Toulouse): Temporary Teaching and Research Attaché (2024–2025). Université de Perpignan Via Domitia: Associate Professor (from 2025). PROMES-CNRS Laboratory (Perpignan, France): Researcher (from 2025).

- **Additional Information:** During his doctoral studies, Dr. Sepúlveda-Oviedo co-developed a patent in partnership with industry, focused on embedded monitoring systems using artificial intelligence for fault diagnosis in photovoltaic installations. In 2023, he was awarded a thesis prize in France recognizing the excellence and impact of his doctoral research. He currently serves as a scientific evaluator for the France 2030 program in areas related to artificial intelligence. He is also actively involved in the international Stic AmSud HAMADI 4.0 project, developing next-generation AI algorithms for diagnosis for Industry 4.0 applications. His research interests include artificial intelligence, predictive maintenance, embedded systems, battery modeling, and renewable energy systems, with a particular focus on autonomous vehicles and distributed energy applications. He has also worked extensively on anomaly detection and modeling of complex dynamic systems.

## Johan Sebastian Suarez Sepulveda

- **Academic Background:** Mechatronics engineer from the National University of Colombia, graduated in 2025, with an emphasis on system modeling and control system design. Also holds a master's degree in Automation and Electronic Control from INSA Toulouse, where the thesis was carried out in collaboration with the Renault Group in France. Currently pursuing a master's degree in industrial Automation in Colombia, with a specialized focus on control systems and process automation.

- **Institutional Affiliations:** Former intern at the Laboratory for Analysis and Architecture of Systems (LAAS-CNRS), Toulouse, France, in the Energy Management Systems Integration Group (ISGE). Currently a master's student in Industrial Automation at the National University of Colombia.

- **Additional Information:** Enthusiastic about dynamic system modeling and control systems. He has developed over 15 academic projects involving prototyped systems for characterization and control, demonstrating strong skills in simulation environments such as Simscape and Simulink. Currently, he continue to conduct research in robust, optimal, and adaptive control systems.

Corinne Alonso

- **Academic Background:** Corinne Alonso is currently Full Professor at the Department of Electrical, Electronic and Automatic Engineering of Toulouse University and has been conducting her research at the Laboratory for Analysis and Architecture of Systems of the French National Centre for Scientific Research (LAAS-CNRS) in Toulouse, France. She was graduated in Power Electronics from Université Paul Sabatier in 1991. Then, she obtained her Ph.D. from Institut National Polytechnique de Toulouse in 1994 and her HdR (Habilitation diriger les Recherches) from the Paul Sabatier University in 2003. She has more than 30 years of academic and professional experience devoted to power electronics and energy conversion, with strong expertise in renewable energy systems and new grids with intermittence energy sources and their impact on resillient micro-grids.

- **Relevant Publications:**

  - A. Gutierrez, Michael Bressan, Fernando Jimenez, Corinne Alonso. Real-time emulation of boost inverter using the Systems Modeling Language and Petri nets. *Mathematics and Computers in Simulation*, 2018, 216-234.

  - Margot Gaetani-Liseo, Corinne Alonso, Bruno Jammes. Identification of ESS Degradations Related to their Uses in Micro-Grids: application to a building lighting network with VRLA batteries. European Journal of Electrical Engineering, 2021, 23 (6), pp.455-466.

  - Antonino Sferlazza, Carolina Albea-Sanchez, Luis Martínez-Salamero, Germain Garcia, Corinne Alonso. Min-Type Control Strategy of a DC-DC Synchronous Boost Converter. *IEEE Transactions on Industrial Electronics*, 2019, 67

- **Institutional Affiliations:** Université Paul Sabatier (Toulouse, France): Graduate in Power Electronics (1991). Institut National Polytechnique de Toulouse: Ph.D. in Power Electronics (1994). Laboratory for Analysis and Architecture of Systems (LAAS-CNRS): Researcher (since 1996). Université Toulouse III – Paul Sabatier: Full Professor (since 2007).

- **Additional Information:** Professor Alonso has led numerous projects in renewable energy integration, with a special focus on photovoltaic systems, converter design, and MPPT control strategies under partial shading and in micro-grids. She was a key contributor to the BIPV Living Lab at LAAS-CNRS, deploying and expanding rooftop PV systems from 100 kWp to 150 kWp in 2022. She co-lead the transverse Energy axis of the LAAS-CNRS and serves on the Board of the international CNRS collaborative project "NextPV" with the RCAST Lab in Japan. Her research activities combine modeling, experimentation, and technology transfer in the fields of smart microgrids, wide-bandgap power devices, and sustainable energy infrastructures (more than 28 PhD supervised, 42 scientific journals and more than 130 publications on international journals, 10 patents).

## Bruno Jammes

- **Academic Background:** Bruno Jammes is Associate Professor at Université de Toulouse. He holds a Ph.D. in Electrical Engineering and has developed extensive expertise in energy systems, particularly in battery modeling, energy storage integration, and microgrid management. His academic background is focused on the design, control, and monitoring of battery-based systems for smart and sustainable energy infrastructures.

- **Relevant Publications:**

  -Jammes, B., Sepúlveda-Oviedo, E.H., Alonso, C. (2025). Suivi temps réel du SoH de batteries lithium-ion. Symposium de Génie Électrique SGE 2025, Toulouse, France.

  -Jammes, B., Anvari-Moghaddam, A.M., Dulout, J., Alonso, C., Guerrero, J.M. (2018). Optimal Design and Operation Management of Battery-Based Energy Storage Systems (BESS) in Microgrids. In: Advancements in Energy Storage Technologies, InTech, pp. 147.

- **Institutional Affiliations:** Université de Toulouse: Associate Professor. Laboratory for Analysis and Architecture of Systems (LAAS-CNRS, Toulouse, France).

- **Additional Information:** His research focuses on the energy management of microgrids with high penetration of renewable sources and energy storage systems. He has contributed to several national and international research projects related to lithium-ion battery aging, SoH (State of Health) monitoring, hybrid energy storage systems, and optimal system design. His work bridges physical modeling, data-driven techniques, and system-level integration to improve reliability and performance of distributed energy infrastructures.